\input amstex

\documentstyle{amsppt}

\hsize450pt

\input epsf

\def\B{\mathop{\text{\rm B}}}
\def\G{\mathop{\text{\rm G}}}
\def\Re{\mathop{\text{\rm Re}}}
\def\tr{\mathop{\text{\rm tr}}}
\def\SU{\mathop{\text{\rm SU}}}
\def\S{\mathop{\text{\rm S}}}
\def\PU{\mathop{\text{\rm PU}}}
\def\Im{\mathop{\text{\rm Im}}}
\def\ta{\mathop{\text{\rm ta}}}
\def\Arg{\mathop{\text{\rm Arg}}}

\hsize450pt

\topmatter\title Complex Hyperbolic Structures on Disc Bundles over
Surfaces.\\ II. Example of a Trivial Bundle\endtitle\author Sasha
Anan$'$in and Nikolay Gusevskii\endauthor
\thanks\endthanks
\thanks Supported by CNPq\endthanks
\address Departamento de Matem\'atica, IMECC, Universidade Estadual de
Campinas,\newline13083-970--Campinas--SP, Brasil\endaddress
\email Ananin$_-$Sasha\@yahoo.com\endemail
\address Departamento de Matem\'atica, ICEX, Universidade Federal de
Minas Gerais,\newline31161-970--Belo Horizonte--MG, Brasil\endaddress
\email Nikolay\@mat.ufmg.br\endemail
\abstract This article is based on the methods developed in [AGG]. We
construct a complex hyperbolic structure on a trivial disc bundle over
a closed orientable surface $\Sigma$ (of genus $2$) thus solving a long
standing problem in complex hyperbolic geometry (see [Gol1, p.~583] and
[Sch, p.~14]). This example answers also [Eli, Open Question 8.1] if a
trivial circle bundle over a closed surface of genus $>1$ admits a
holomorphically fillable contact structure. The constructed example $M$
satisfies the relation $2(\chi+e)=3\tau$ which is necessary for the
existence of a holomorphic section of the bundle, where
$\chi=\chi\Sigma$ stands for the Euler characteristic of $\Sigma$,
$e=eM$, for the Euler number of the bundle, and $\tau$, for the Toledo
invariant. (The relation is also valid for the series of examples
constructed in [AGG].) Open question: Does there exist a holomorphic
section of the bundle $M$?
\endabstract\endtopmatter\document

\rightheadtext{II. Example of a Trivial Bundle}

\bigskip

\centerline{\bf1. Introduction}

\bigskip

In this article, we answer a long standing question raised by several
authors (for instance, [Eli, Open Question 8.1], [Gol1, p.~583], and
[Sch, p.~14]): Does there exist a complex hyperbolic structure on a
trivial disc bundle $M$ over a closed orientable surface $\Sigma$? In
what follows, we construct such an example.

Our start point was the relation $2(\chi+e)=3\tau$ valid for the series
of examples constructed in [AGG], where $\chi=\chi\Sigma$ stands for
the Euler characteristic of $\Sigma$, $e=eM$, for the Euler number of
the bundle, and $\tau$, for the Toledo invariant of the representation
$\pi_1\Sigma\to\PU(2,1)$ defined by $M$ ($\chi$, $e$, and $\tau$ are
negative). It is known that $\tau\in\frac{2}{3}\Bbb Z$ and
$|\tau|\le|\chi|$ [Tol]. If the relation $2(\chi+e)=3\tau$ is valid,
$\Sigma$ is of genus $2$, and $\tau$ is not integer, then $e=0$
automatically. Even if $2(\chi+e)\ne3\tau$, an example of a complex
hyperbolic disc bundle $M$ over $\Sigma$ of genus $2$ and noninteger
$\tau M$ is curious by itself since there is another challenging
conjecture, a complex hyperbolic variant of the question raised in
[GLT]: Is the inequality $|eM|\le|\chi\Sigma|$ valid for every
(complex) hyperbolic disc bundle $M$ over a closed orientable surface
$\Sigma$?

In the complex hyperbolic plane $\Bbb H_\Bbb C^2$, we find three
reflections $R_1,R_2,R_3$ --- the first two are reflections in points,
the third one is the reflection in a real plane --- subject to the
relation $(R_3R_1R_2R_3R_2R_1)^2=1$. A suitable subgroup of index $8$
in the discrete group $T$ generated by $R_1,R_2,R_3$ provides a disc
bundle $M$ in question. In order to prove the discreteness of $T$, we
use the methods developed in [AGG] and show that a suitable elliptic
transversal triangle of bisectors with common slices (see [AGG]) yields
a fundamental polyhedron for $T$. The Toledo invariant of $M$ is not
integer since $(R_3R_1R_2R_3R_2R_1)^2\ne1$ in terms of $\SU(2,1)$. As
expected, the example satisfies the relation $2(\chi+e)=3\tau$. This
makes more intriguing the question raised in [AGG] about the existence
of a holomorphic section.

\medskip

{\bf Acknowledgements.} The authors are very grateful to Carlos
Henrique Grossi Ferreira for his stimulating interest to this work and
for many fruitful conversations on the subject. We dedicate this
article to the memory of Andrei Tyurin.

\bigskip

\centerline{\bf2. The Triangle and the Cake}

\bigskip

In this section, we construct a complex hyperbolic disc bundle over an
orientable surface of genus $3$. In the sequel, the roman I indicates a
reference to [AGG].

\medskip

{\bf2.1. Preliminaries: Notation and Conventions.} We remind some of
the concepts and notation from [AGG] where further details can be
found.

\smallskip

We denote by $V$ a three-dimensional $\Bbb C$-vector space equipped
with a hermitian form $\langle-,-\rangle$ of signature $++-$. The open
$4$-ball and the $3$-sphere
$$\B V:=\big\{p\in\Bbb{CP}V\mid\langle p,p\rangle<0\big\},\qquad\S
V:=\big\{p\in\Bbb{CP}V\mid\langle p,p\rangle=0\big\}$$
are the {\it complex hyperbolic plane\/} $\Bbb H_\Bbb C^2$ and its
ideal boundary, the {\it absolute\/} $\partial\Bbb H_\Bbb C^2$. We put
$\overline{\B}V:=\B V\cup\S V$. Every projective line $L$ in
$\Bbb{CP}V$ has the form $\Bbb{CP}p^\perp$, where
$p^\perp:=\big\{v\in V\mid\langle v,p\rangle=0\big\}$. We say that $p$
is the {\it polar point\/} to $L$. If $p\notin\overline\B V$, then
$\Bbb{CP}p^\perp\cap\overline\B V$ is a {\it complex geodesic.}

The reflection
$$R(p):x\mapsto2\frac{\langle x,p\rangle}{\langle p,p\rangle}p-x$$
is defined for every nonisotropic $p\in\Bbb{CP}V$ and belongs to
$\SU V$. If $p\in\B V$, then $R(p)$ is the reflection in $p$. If
$p\notin\overline\B V$, then $R(p)$ is the reflection in the complex
geodesic $\Bbb{CP}p^\perp\cap\overline\B V$.

For $b,e\in\overline\B V$, we denote by $\G[b,e]$ the closed segment of
geodesic oriented from $b$ to $e$ and by $\G{\prec}b,e{\succ}$, the
corresponding oriented geodesic in $\overline\B V$. (Thus, the geodesic
$\G{\prec}b,e{\succ}$ includes its vertices.) Similar notation applies
to bisectors. Given two ultraparallel complex geodesics $C_1$ and
$C_2$, they determine a unique oriented segment of bisector
$B[C_1,C_2]$ which begins with the slice $C_1$ and ends with the
slice~$C_2$. The full oriented bisector related to the segment
$B[C_1,C_2]$ is denoted by $B{\prec}C_1,C_2{\succ}$. In the same
manner, we define a semifull bisector $B[C_1,C_2{\succ}$.

Let $B$ and $B'$ be two oriented bisectors with a common slice $S$ and
let $p\in S$. Let $n$ and $n'$ stand for the normal vectors at $p$ to
$B$ and to $B'$, respectively. Since both $n$ and $n'$ are tangent to
the naturally oriented complex geodesic passing through $p$ and
orthogonal to $S$, it makes sense to measure the oriented angle from
$n$ to $n'$. This angle is said to be the {\it oriented angle\/} from
$B$ to $B'$ at $p$.

Three pairwise ultraparallel complex geodesics $C_1,C_2,C_3$ give raise
to the {\it triangle of bisectors\/} $\Delta(C_1,C_2,C_3)$ whose {\it
vertices\/} are  $C_1,C_2,C_3$ and whose {\it sides\/} are
$B[C_1,C_2]$, $B[C_2,C_3]$, $B[C_3,C_1]$. Let $\Delta(C_1,C_2,C_3)$ be
a transversal triangle oriented in the counterclockwise sense (see
Subsection I.2.5). By~Lemma I.2.2.2, the full bisectors
$B{\prec}C_1,C_2{\succ}$ and $B{\prec}C_1,C_3{\succ}$ divide
$\overline\B V$ into four {\it sectors.} The sector that includes
$\Delta(C_1,C_2,C_3)$ is said to be the {\it interior sector\/} (or
simply the {\it sector\/}) of $\Delta(C_1,C_2,C_3)$ at $C_1$. This
sector is {\it oriented\/} from its {\it first side\/}
$B[C_1,C_2{\succ}$ to its {\it second side\/} $B[C_1,C_3{\succ}$. By
this definition, the sector includes its sides. The corresponding
sides $B[C_1,C_2]$ and $B[C_1,C_3]$ of the triangle are respectively
included in the sides of the sector at $C_1$.

\medskip

{\bf2.2. The Triangle.} In this subsection, we describe a triangle of
bisectors $\Delta$ which is a building block for our example. We just
list the elements and the properties of the triangle and postpone until
Sections 4 and 5 the proof of the existence.

\smallskip

$\bullet$ Three pairwise ultraparallel complex geodesics $C_1,C_2,C_3$
provide the elliptic transversal triangle of bisectors
$\Delta:=\Delta(C_1,C_2,C_3)$ oriented in the counterclockwise sense.
In terms of the polar points $p_i\notin\overline\B V$, we have
$C_i=\Bbb{CP}p_i^\perp\cap\overline\B V$.

$\bullet$ Let $\G[b_i,e_{i+1}]$, $m_i\in\G[b_i,e_{i+1}]$, and
$q_i\notin\overline\B V$ stand respectively for the real spine, its
middle point, and for the polar point to the middle slice $M_i$ of the
segment $B[C_i,C_{i+1}]$ (the indices are modulo $3$). The points
$b_2,e_2\in C_2$ are distinct.

$\bullet$ We denote by $R_0:=R(p_1)$, $R_1:=R(m_1)$, and $R_2:=R(m_2)$
respectively the reflections in $C_1$, in~$m_1$, and in $m_2$. Let
$\vartheta:=\exp\displaystyle\frac{\pi i}3$. The points
$$c_1:=e_1,\qquad c_2:=R_1c_1,\qquad c_3:=R_2c_2,\qquad d_3:=b_3,\qquad
d_2:=R_2d_3,\qquad d_1:=R_1d_2$$
satisfy the relation
$\displaystyle\frac{\langle c_1,c_3\rangle\langle
c_3,c_2\rangle}{\langle c_1,c_2\rangle}\in\Bbb R\overline\vartheta i$.
We define the oriented geodesic segments
$$S_1:=\G[c_1,c_2],\ S_2:=\G[c_2,c_3],\ S_3:=\G[c_3,c_1],\
T_1:=\G[d_1,d_2],\ T_2:=\G[d_2,d_3],\ T_3:=\G[d_3,d_1]$$

\vskip5pt

\noindent
$\hskip243pt\vcenter{\hbox{\epsfbox{Pict.1}}}$

\rightskip220pt

\vskip-178pt

\noindent
that form two closed\footnote{For two curves $\gamma_1$ and $\gamma_2$
such that $\gamma_2$ begins with the end of $\gamma_1$, we denote by
$\gamma_1\cup \gamma_2$ their path-product.}
curves $\sigma:=S_1\cup S_2\cup S_3$ and
$\varsigma:=T_1\cup T_2\cup T_3$ in $\Delta$. The curves $\sigma$ and
$\varsigma$ have at most 3 common points.

$\bullet$ Let $\beta_i\in[0,2\pi)$ denote the oriented angle from
$B[C_i,C_{i+1}]$ to $B[C_i,C_{i-1}]$ at $c_i$ (the indices are modulo
$3$). These angles satisfy the relation
$0<\beta_1+\beta_2+\beta_3<\pi$.

$\bullet$ There exist a meridian $\Gamma$ of $B[C_3,C_1]$ and a vertex
$w_3\in\S V$ of the geodesic $\Gamma\cap C_3$ such that
$c_3,d_1\in\Gamma$ and $R(m_3)w_3\ne R_1R_2w_3$. The point
$R(q_1)R(q_3)w_3$ lies in $C_2$ on the side of the normal vector to the
oriented geodesic $\G{\prec}b_2,e_2{\succ}$. Obviously,
$S_3,T_3\subset\Gamma$.

$\bullet$ The relation $R_3R_1R_2R_3R_2R_1R_0=\vartheta^{-2}$ holds in
$\widehat\SU V$, where $R_3$ stands for the reflection in $\Gamma$ and
$\widehat\SU V$\penalty-10000

\vskip-12pt

\rightskip0pt

\noindent
denotes the group $\SU V$ extended with a reflection in a real plane.

$\bullet$ The triangles $\Delta$ and $R_3\Delta$ lie on different sides
from the full bisector $B{\prec}C_3,C_1{\succ}$.

\smallskip

By Theorem I.2.5.1, $\Delta$ and $\S V$ bound a fibred polyhedron
$P\subset\overline\B V$ which lies on the side of the normal vector to
each $B[C_i,C_{i+1}]$, $i=1,2,3$. Clearly,
$\Delta\cap\overline\B V=\partial_0P$. (See Subsection I.2.4 for
definitions.)

\medskip

{\bf2.3. The Cake.} The example we are going to construct mimics the
following example in the hyperbolic plane $\Bbb H_\Bbb R^2$.

\smallskip

Let $T$ be a geodesic triangle in $\Bbb H_\Bbb R^2$ with vertices
$v_1,v_2,v_3$ listed in the counterclockwise order and with interior
angles $\alpha_1,\alpha_2,\alpha_3$ satisfying the relation
$\alpha_1+\alpha_2+\alpha_3=\frac\pi2$. The edge $s_i$ of $T$ begins
with $v_i$ and ends with $v_{i+1}$ (the indices are modulo $3$).
We choose the clockwise orientation of the triangle~$T$. Hence,
$\partial T=(s_1\cup s_2\cup s_3)'$, where $\xi'$ denotes the curve
$\xi$ taken with the opposite orientation. Let $r_0$, $r_i$, and $r_3$
denote the reflections in $v_1$, in the middle point of $s_i$, and in
$s_3$, respectively, $i=1,2$. By~Poincar\'e's Polyhedron Theorem, $T$
is a fundamental domain for the group $\Upsilon$ generated by the
$r_i$'s, and the defining relations are $r_3r_1r_2r_3r_2r_1r_0=1$ and
$r_i^2=1$, $i=0,1,2,3$. This fact becomes more evident if we throw a
look at the following cake:

\newpage

\noindent
$\hskip0pt\vcenter
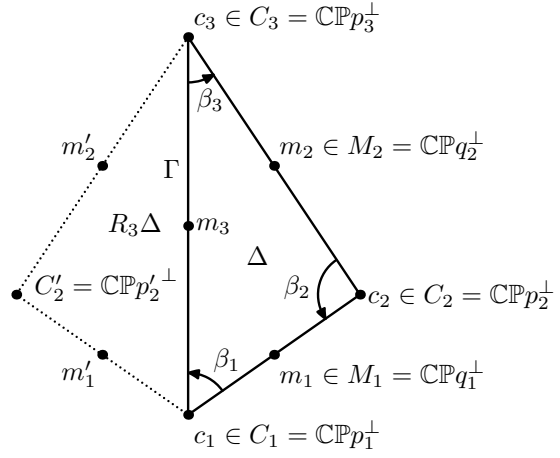
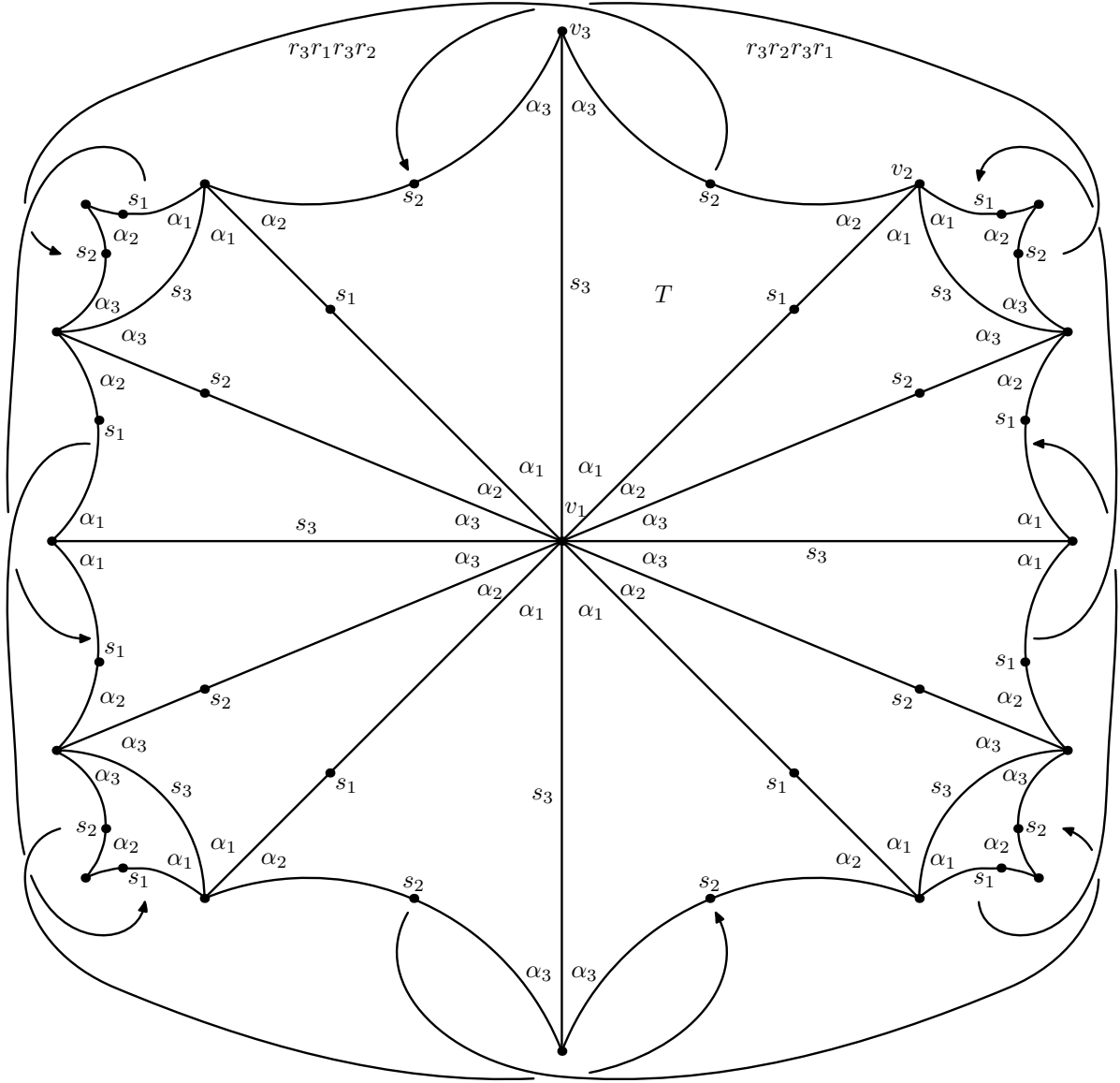
{\hbox{\epsfbox{Pict.2}}}$

\vskip10pt

\noindent
In each triangle, we indicated the angles $\alpha_1,\alpha_2,\alpha_3$
and the type of the edge $s_1,s_2,s_3$ assuming that adjacent triangles
with common $s_i$ are glued with the help of a suitable conjugate to
$r_i$, $i=1,2,3$. The arrows indicate all the identifications of the
edges on the boundary of the cake. Some of them are marked with a
corresponding element of $\Upsilon$. On its boundary, the cake contains
3 cycles of vertices and 8 pairs of edges to be identified. Therefore,
the cake provides a surface $\Sigma$ whose fundamental group
$\pi_1\Sigma\subset\Upsilon$ is of index $16$. The genus of $\Sigma$
equals $3$ since $\chi=3-8+1=2-2g$.

\medskip

{\bf2.4. Toledo invariant.} We define a representation
$\varrho:\Upsilon\to\widehat\PU V$ by $r_i\mapsto R_i$, $i=0,1,2,3$,
which induces a representation
$\varrho|_{\pi_1\Sigma}:\pi_1\Sigma\to\PU V$, where the group
$\widehat\PU V$ stands for $\PU V$ extended with a reflection in a real
plane. Moreover, there exists a $\varrho$-equivariant continuous
mapping $\varphi:\Bbb H_\Bbb R^2\to\B V$ that sends $v_i$ to $c_i$ and
$s_i$ onto $S_i$ since we can map $T$ onto an arbitrary disc
$B'\subset\B V$ with $\partial B'=\sigma'$. (We can even assume that
$B'\subset P$ because $P$ is fibred and, therefore, is a closed
$4$-ball.) It follows a sort of a copy of Proposition I.2.1.6.

\medskip

{\bf2.4.1. Proposition.} {\sl The Toledo invariant of\/
$\varrho|_{\pi_1\Sigma}$ equals\/ $-\frac{8}{3}$.}

\medskip

{\bf Proof.} The Toledo invariant of $\varrho|_{\pi_1\Sigma}$ is
defined as $\tau:=4\displaystyle\frac1{2\pi}\int_C\varphi^*\omega$
[Tol], where $C$ stands for the cake oriented in the clockwise sense
and $\omega$, for the K\"ahler form described in explicit terms in
Lemma~I.2.1.3. Taking $u:=c_1\in\B V$, by Lemma I.2.1.3, we obtain
$$\tau=\frac{32}\pi\int_{B'}\omega=\frac{32}\pi\int_{\sigma'}
P_{c_1}=-\frac{32}\pi\sum\limits_{i=1}^3\int_{S_i}P_{c_1}.$$
By Lemma I.2.1.4, the terms with $i\ne2$ vanish. By Lemmas I.2.1.4 and
I.2.1.3,
$$\int_{S_2}P_{c_1}=\int_{S_2}(P_{c_1}-P_{c_2})=
\int_{S_2}df_{c_1,c_2}.$$
This number is the total variation of
$\displaystyle\frac12\Arg\frac{\langle c_1,x\rangle\langle
x,c_2\rangle}{\langle c_1,c_2\rangle}$,
while $x$ runs over $S_2\subset\B V$ from $c_2$ to $c_3$. By~Lemma
I.2.1.5,
$\displaystyle\frac{\langle c_1,x\rangle\langle x,c_2\rangle}{\langle
c_1,c_2\rangle}$
is never real nonnegative. It follows that
$$\int_{S_2}df_{c_1,c_2}=\frac12\Arg\frac{\langle c_1,c_3\rangle\langle
c_3,c_2\rangle}{\langle c_1,c_2\rangle}-\frac12\Arg\frac{\langle
c_1,c_2\rangle\langle c_2,c_2\rangle}{\langle
c_1,c_2\rangle}=\frac12\Arg\frac{\langle c_1,c_3\rangle\langle
c_3,c_2\rangle}{\langle c_1,c_2\rangle}-\frac\pi2$$
since $\langle c_2,c_2\rangle<0$. From
$0\ne\displaystyle\frac{\langle c_1,c_3\rangle\langle
c_3,c_2\rangle}{\langle c_1,c_2\rangle}\in\Bbb R\overline\vartheta i$
and $\Arg\overline\vartheta i=\frac\pi6$, we obtain two possible values
of $\tau$~: $-\frac83$ and $\frac{40}3$. The latter is impossible
because of $|\tau|\le|\chi|=4$ (see [Tol])
$_\blacksquare$

\medskip

{\bf2.5. Poincar\'e's Polyhedron Theorem.} In terms of $\widehat\PU V$,
the relation between the $R_i$'s can be written as
$R_3R_1R_2R_3R_2R_1=R_0$. Hence,
$$R_3R_1R_2R_3R_2R_1R_3R_1R_2R_3R_2R_1=1.$$
Let $W_i$ denote the subword of length $i$ of the last relation,
$0\le i\le12$. For instance, $W_1:=R_3$, $W_2:=R_3R_1$,
$W_5:=R_3R_1R_2R_3R_2$, etc. Denoting by $\nabla'$ the triangle of
bisectors $\nabla$ taken with opposite orientation, we define the
triangles $\Delta_i$, $1\le i\le12$, as follows:
$$\Delta_i:=W_i\Delta'\text{ \rm for
}i=1,2,3,7,8,9;\qquad\Delta_i:=W_i\Delta\text{ \rm for
}i=4,5,6,10,11,12.$$
Since $R_1$ and $R_2$ preserve the property `to be oriented in the
counterclockwise sense' and $R_3$ alters it, all 12 triangles are
transversal and oriented in the counterclockwise sense. It follows from
$$R_1C_1=C_2,\quad R_1C_2=C_1,\quad R_2C_2=C_3,\quad R_2C_3=C_2,\quad
R_3C_1=C_1,\quad R_3C_3=C_3$$
that
$$W_iC_1=C_1\text{ \rm for }i=1,6,7,12;\qquad W_iC_2=C_1\text{ \rm for
}i=2,5,8,11;\qquad W_iC_3=C_1\text{ \rm for }i=3,4,9,10.$$
Similarly,
$$R_1c_1=c_2,\quad R_1c_2=c_1,\quad R_2c_2=c_3,\quad R_2c_3=c_2,\quad
R_3c_1=c_1,\quad R_3c_3=c_3$$
imply
$$W_ic_1=c_1\text{ \rm for }i=1,6,7,12;\qquad W_ic_2=c_1\text{ \rm for
}i=2,5,8,11;\qquad W_ic_3=c_1\text{ \rm for }i=3,4,9,10.$$

The first and second sides of the sector of $\Delta$ at $C_1$ equal
respectively the second side of the sector of $R_1\Delta$ at $C_1$ and
the first side of the sector of $R_3\Delta'$ at $C_1$. The first and
second sides of the sector of $\Delta$ at $C_2$ equal respectively the
second side of the sector of $R_2\Delta$ at $C_2$ and the first side of
the sector of $R_1\Delta$ at $C_2$. The first and second sides of the
sector of $\Delta$ at $C_3$ equal respectively the second side of the
sector of $R_3\Delta'$ at $C_3$ and the first side of the sector of
$R_2\Delta$ at $C_3$.

These facts immediately imply that the second side of the sector of
$\Delta_{i-1}$ at $C_1$ coincides with the first side of the sector of
$\Delta_i$ at $C_1$ for every $0\le i\le11$
($\Delta_0:=\Delta_{12}=\Delta$) --- we denote this common side of the
sectors by $L_i$.

From the relation $R_3R_1R_2R_3R_2R_1=R_0$ and from the fact that,
reflecting a semifull bisector in its initial slice, we obtain the
other part of the full bisector, we deduce that $R_0L_i=L_{i+6}$ and
that $L_i$ and $L_{i+6}$ form a full bisector $B_i$, $1\le i\le6$.

By the above, $B_i$ and $B_{i+1}$ (the indices are modulo $6$) are
transversal along their common slice $C_1$ and the angle from $L_i$ to
$L_{i+1}$ (the indices are modulo $12$) at $c_1\in C_1\cap\B V$ equals
$$\beta_1\text{ \rm for }i=1,6,7,12;\qquad\beta_2\text{ \rm for
}i=2,5,8,11;\qquad\beta_3\text{ \rm for }i=3,4,9,10.$$

{\bf2.5.1. Lemma.} {\sl$\beta_1+\beta_2+\beta_3=\frac\pi2$. The
bisectors\/ $B_i$'s are pairwise transversal along their common slice\/
$C_1$.}

\medskip

{\bf Proof} explores the idea of Lemma I.2.3.1. (The picture close to
that lemma is also useful.) Let $C$ be a complex geodesic orthogonal to
$C_1$ and passing through $c_1$. Each bisector $B_i$ intersects it in a
smooth curve divided by $c_1$ into two parts $l_i$ and $l_{i+6}$. The
oriented angle at $c_1$ from $l_i$ to $l_{i+1}$ (the indices are modulo
$12$) is the one indicated above Lemma 2.5.1. This implies that
$4\beta_1+4\beta_2+4\beta_3=2n\pi$, $n\in\Bbb N$. Since
$0<\beta_1+\beta_2+\beta_3<\pi$, we obtain the first assertion. Now, we
can scan our configuration with a complex geodesic $C$ orthogonal to
$C_1$ making it pass through other points in $C_1\cap\B V$. The sum of
the angles from $l_i$ to $l_{i+1}$ has to be the same, $2\pi$, and
either of the angles takes values in $(\varepsilon,\pi)$ for a suitable
$\varepsilon>0$ since $B_i$ and $B_{i+1}$ are transversal. This implies
that the angle from $l_i$ to $l_j$ exceeds $\varepsilon$ unless
$i\equiv j\mod12$
$_\blacksquare$

\medskip

We have proven more than claimed in Lemma 2.5.1. The triangles
$\Delta_i$'s are separated by the sides $L_i$'s of their sectors at
$C_1$. Looking at the cake, we can see 4 more triangles: the triangle
$\Delta_{13}:=W_2R_3\Delta$ having a common side with $\Delta_2$, the
triangle $\Delta_{14}:=W_5R_3\Delta'$ having a common side with
$\Delta_5$, the triangle $\Delta_{15}:=W_8R_3\Delta$ having a common
side with $\Delta_8$, and the triangle $\Delta_{16}:=W_{11}R_3\Delta'$
having a common side with $\Delta_{11}$. It is easy to explicitly find
an isometry mapping the pair of adjacent triangles $\Delta_{j+13}$ and
$\Delta_{3j+2}$ onto a suitable pair $\Delta_i$ and $\Delta_{i+1}$ so
that an arbitrarily chosen vertex of the common side of the first pair
goes to $C_1$, $j=0,1,2,3$. This allows us to conclude that the
triangle $\Delta_{j+13}$ is included in the sector of $\Delta_{3j+2}$
at $C_1$. Therefore, the sides of the triangles $\Delta_i$,
$1\le i\le16$, that remain `unglued' form a simple configuration
$\Cal K$ of bisectors (the orientation of $\Cal K$ comes from the
triangles $\Delta_i$ oriented in the counterclockwise sense) such that
every two adjacent bisectors in $\Cal K$ are transversal along their
common slice. (See Subsection I.2.3 for definitions.)

Moreover, the total angle of every cycle of vertices of $\Cal K$ equals
$2\pi$. To measure the total angle of a cycle of vertices, we take a
point $c_0\in C_0\cap\B V$, where $C_0$ is some vertex of the cycle of
vertices, and, with the help of the isometries involved in the cycle,
obtain from $c_0$ a point $c\in C\cap\B V$ in each vertex $C$ of the
cycle. We define the total angle of the cycle of vertices as the sum of
oriented angles at the $c$'s between corresponding bisectors. The
reasons similar to those used in the proof of Lemma 2.5.1 (or~Lemma
I.2.3.1) show that if the total angle equals $2\pi$ for some $c_0$, it
is independent of the choice of~$c_0$.

\smallskip

By Theorem I.2.5.1, each $\Delta_i$ together with $\S V$ bounds a
fibred polyhedron $P_i$ located on the side of the normal vectors to
the sides of $\Delta_i$. By Remark I.2.4.3, $\Cal K$ and $\S V$ bound a
fibred polyhedron $\Cal C$, called 4-{\it Cake,} that lies on the side
of the normal vectors to each bisector involved in $\Cal K$. We denote
by $G$ the group generated by the 8 isometries that identify 16
bisectors on the boundary $\partial_0\Cal C=\Cal K\cap\overline\B V$ :

\smallskip

\noindent
$I_1:=W_2R_3R_2W_0^{-1}$ maps the side $B[C_2,C_3]$ of
$\Delta_{12}=\Delta$ onto the side $B[W_2R_3C_3,W_2R_3C_2]$ of
$\Delta_{13}$,

\noindent
$I_2:=W_4R_1R_3W_2^{-1}$ maps the side $B[W_2R_3C_1,W_2R_3C_2]$ of
$\Delta_{13}$ onto the side $B[W_4C_2,W_4C_1]$ of $\Delta_4$,

\noindent
$I_3:=W_5R_3R_1W_3^{-1}$ maps the side $B[W_3C_2,W_3C_1]$ of $\Delta_3$
onto the side $B[W_5R_3C_1,W_5R_3C_2]$ of $\Delta_{14}$,

\noindent
$I_4:=W_7R_2R_3W_5^{-1}$ maps the side $B[W_5R_3C_3,W_5R_3C_2]$ of
$\Delta_{14}$ onto the side $B[W_7C_2,W_7C_3]$ of $\Delta_7$,

\noindent
$I_5:=W_8R_3R_2W_6^{-1}$ maps the side $B[W_6C_2,W_6C_3]$ of $\Delta_6$
onto the side $B[W_8R_3C_3,W_8R_3C_2]$ of $\Delta_{15}$,

\noindent
$I_6:=W_{10}R_1R_3W_8^{-1}$ maps the side $B[W_8R_3C_1,W_8R_3C_2]$ of
$\Delta_{15}$ onto the side $B[W_{10}C_2,W_{10}C_1]$ of $\Delta_{10}$,

\noindent
$I_7:=W_{11}R_3R_1W_9^{-1}$ maps the side $B[W_9C_2,W_9C_1]$ of
$\Delta_9$ onto the side $B[W_{11}R_3C_1,W_{11}R_3C_2]$ of
$\Delta_{16}$,

\noindent
$I_8:=W_1R_2R_3W_{11}^{-1}$ maps the side
$B[W_{11}R_3C_3,W_{11}R_3C_2]$ of $\Delta_{16}$ onto the side
$B[W_1C_2,W_1C_3]$ of $\Delta_1$.

\smallskip

Now, we can read almost literally the proofs of Remark 4.1 and Lemmas
I.4.2, I.4.4, and I.4.7 adapting them to our case. The analogs of
Lemmas I.2.3.1 and I.4.8 are already proven and we can simply repeat
the proof of Theorem I.2.2.3 adapted to our situation.

By Lemma I.2.4.1, we can extend the slice bundle of the boundary
$\partial_0\Cal C$ to some disc bundle of $\Cal C$. Since the
generators of $G$ preserve the slice bundle of $\partial_0\Cal C$, we
arrive at the

\medskip

{\bf2.5.2. Proposition.} {\sl The\/ $4$-Cake $\Cal C$ is a fibred
polyhedron that provides a fundamental polyhedron for the group\/
$G\subset\PU V$ isomorphic to\/ $\pi_1\Sigma$, where\/ $\Sigma$ stands
for a closed orientable surface of genus\/ $3$. The complex hyperbolic
manifold\/ $M:=\B V/G$ is a disc bundle over\/ $\Sigma$
$_\blacksquare$}

\medskip

A couple of remarks about Poincar\'e's Polyhedron Theorem:

1. To prove a $4$-dimensional version of Poincar\'e's Polyhedron
Theorem, even as plane as ours, it is essential to have some estimate
concerning the behaviour of faces that is uniform along the common part
of the faces. (In our case, such estimate is provided by Lemma I.4.3.)
Otherwise, it is easy to construct a counterexample.

2. Clearly, there is a general analog of plane Poincar\'e's Polyhedron
Theorem that allows to prove the discreteness of the group generated by
$R_1,R_2,R_3$ without gluing the $4$-Cake $\Cal C$ from the polyhedra
$P_i$'s related to the triangles $\Delta_i$'s (see [AGr]).

\bigskip

\centerline{\bf3. A Piece of Elementary Topology and a Piece of the
Cake}

\bigskip

In this section, we find the Euler number $eM$ of the bundle $M$
constructed in Proposition 2.5.2.

\smallskip

Denote by $E_i$ (by $F_i$) the complex geodesic spanned by the segment
of geodesic $S_i$ (by the segment of geodesic $T_i$), $i=1,2$. Since
each of $S_i$ and $T_i$ is not included into a slice of
$B[C_i,C_{i+1}]$, these complex geodesics are transversal to
$B[C_i,C_{i+1}]$. The segments $S_i$ and $T_i$ intersect only in $m_i$.
Therefore, $E_i$~and $F_i$ are distinct and, therefore, transversal.

The closed curve
$\sigma=S_1\cup S_2\cup S_3\subset\partial_0P=\Delta\cap\overline\B V$
is contractible in $P$. Let $B\subset P$ denote a closed simple disc
such that $\partial B=\sigma$. We orient $B$ with respect to the
orientation of $\sigma$. We can assume that, in some open neighbourhood
of $m_i$, the disc $B$ coincides with $P\cap E_i$, $i=1,2$. Since
$R_iS_i=S_i$ for $i=1,2,3$, some number of discs conjugate to $B$
(altering the orientation each time we apply $R_3$), glued along the
geodesics on their boundaries that are conjugate to $S_i$'s, produces a
section $\Cal B$ of the disc bundle $M$. By construction, $\Cal B$
coincides with a suitable conjugate to $E_i$ in an open neighbourhood
of a conjugate to $m_i$, $i=1,2$.

In a similar way, we can construct one more section of the bundle
starting with the closed curve
$\varsigma=T_1\cup T_2\cup T_3\subset\partial_0P=\Delta\cap\overline\B
V$
contractible in $P$. As above, we can find a closed simple disc
$D\subset P$ such that $\partial D=\varsigma$ and orient $D$ with
respect to the orientation of $\varsigma$. We can assume that $D$
coincides with $P\cap F_i$ in some open neighbourhood of $m_i$,
$i=1,2$. Moreover, the fact that $\sigma$ and $\varsigma$ intersect in
a finite number of points imply that we can take $D$ transversal to
$B$. The disc $D$ generates a section $\Cal D$ of the disc bundle $M$
which coincides with a suitable conjugate to $F_i$ in an open
neighbourhood of a conjugate to $m_i$, $i=1,2$.

By construction, $\Cal B$ and $\Cal D$ are transversal. All we need in
order to find $eM$ is to count $\#\Cal B\cap\Cal D$, the intersection
number of $\Cal B$ and $\Cal D$.

It is convenient to count first the intersection number that
corresponds to the intersections lying in the {\it piece\/}
$\Cal P:=P\cup R_3P$ of the $4$-Cake. It is easy to see that $R_3$
preserves the orientation of $\Bbb{CP}V$ and alters the natural
orientation of every complex geodesic. The reflections $R_1$ and $R_2$
preserve the orientation of $\Bbb{CP}V$ and the natural orientation of
every complex geodesic. We denote by $n$ the intersection number of
$B\cup R_3B'$ and $D\cup R_3D'$ in the interior of $\Cal P$, where $X'$
denotes the disc $X$ taken with opposite orientation. By construction,
$\Cal B$ and $\Cal D$ really intersect in $m_i$, $i=1,2$; we denote the
corresponding number by $\varepsilon_i$. Obviously,
$\varepsilon_i=\pm1$ depending on two relations. One is between the
orientation of $\Cal B$ and the natural orientation of $E_i$, the other
is between the orientation of $\Cal D$ and the natural orientation
of~$F_i$. If both coincide or both differ, we have $\varepsilon_i=1$.
Otherwise, $\varepsilon_i=-1$. The same relations are valid for the
orientations induced by $R_3$ in the pairs $(R_3\Cal B,R_3E_i)$ and
$(R_3\Cal D,R_3F_i)$. The relations will be the same if we
simultaneously alter them all. We can conclude that the intersection
number of $\Cal B$ and $\Cal D$ at $m'_i:=R_3m_i$ equals
$\varepsilon_i$.

\medskip

{\bf3.1. Proposition.} {\sl$eM\equiv0\mod8$.}

\medskip

{\bf Proof.} The 8 polyhedra conjugate to ${\Cal P}$ tessellate
$\Cal C$. Hence, $eM=8(n+\varepsilon_1+\varepsilon_2)$
$_\blacksquare$

\medskip

If $eM\ne0$, then we have constructed an example solving the complex
hyperbolic variant of the conjecture in [GLT] which asks if the
inequality $|eM|\le|\chi\Sigma|$ is valid for every hyperbolic disc
bundle $M$ over a closed orientable surface $\Sigma$.

\smallskip

In order to find the exact value of $eM$, we examine the pieces of the
$4$-Cake more closely and construct one more section.

\smallskip

We define
$$w_2:=R_2w_3\in C_2\cap\S V,\qquad w_1:=R_1w_2\in C_1\cap\S V.$$
It follows from $R_0w_1=w_1$, $R_0=R_3R_1R_2R_3R_2R_1$, and
$R_3w_3=w_3$ that
$$w_1=R_0w_1=R_3R_1R_2R_3R_2R_1w_1=R_3R_1R_2R_3w_3=R_3R_1R_2w_3=
R_3w_1,$$
i.e., $w_1\in\Gamma$.

\noindent
$\hskip263pt\vcenter{\hbox{\epsfbox{Pict.3}}}$

\rightskip200pt

\vskip-150pt

We claim that $R(q_3)w_3=w_1$. Every meridian of a bisector intersects
the boundary of every slice in two points. Also, every meridian of a
bisector is stable under the reflections in slices and under the
reflections in points in the real spine of the bisector. Therefore,
both $R(m_3)$ and $R(q_3)$ map the two-point set
$\Gamma\cap C_3\cap\S V$ onto $\Gamma\cap C_1\cap\S V$. Clearly, these
two mappings are different.\footnote{$R(m_3)R(q_3)$ is the reflection
in the complex spine of $B[C_3,C_1]$ and, therefore, induces in
$\Gamma$ the reflection in the real spine of $B[C_3,C_1]$.}
Since $R(m_3)w_3\ne R_1R_2w_3=w_1$, we obtain $R(q_3)w_3=w_1$.

The cylinder $B[C_3,C_1]\cap\S V$ and the meridian $\Gamma$ intersect
in two curves (ideal meridional curves in terms of the beginning of
Subsection I.2.1). By Lemma I.2.1.1, we get a curve
$\ell\subset\Gamma\cap B[C_3,C_1]\cap\S V$ that begins with $w_3$ and
ends with $w_1$.

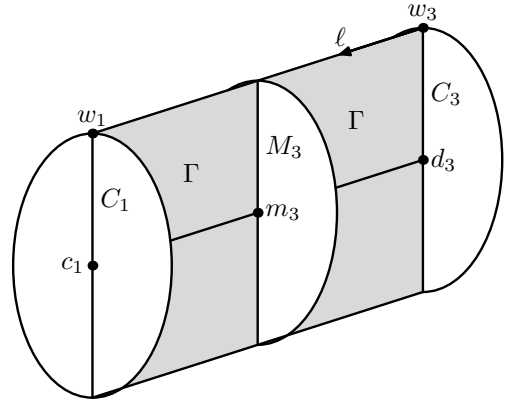
\rightskip0pt

We will denote by $\xi'$ the curve $\xi$ (not necessarily closed) taken
with the opposite orientation.

Let $\mu_1\subset B[C_1,C_2]\cap\S V$ be a simple curve that begins
with $w_1$, ends with $w_2$, and is transversal to both $C_1\cap\S V$
and $C_2\cap\S V$. Let $\nu_1\subset B[C_2,C_3]\cap\S V$ be a simple
curve that begins with $w_2$, ends with $w_3$, and is transversal to
both $C_2\cap\S V$ and $C_3\cap\S V$. Obviously, the curve
$\mu_2:=R_1\mu'_1\subset B[C_1,C_2]\cap\S V$ begins with $w_1$, ends
with $w_2$, and is transversal to both $C_1\cap\S V$ and $C_2\cap\S V$.
The curve $\nu_2:=R_2\nu'_1\subset B[C_2,C_3]\cap\S V$ begins with
$w_2$, ends with $w_3$, and is transversal to both $C_2\cap\S V$ and
$C_3\cap\S V$. We~obtain two closed curves
$\gamma_i:=\mu_i\cup\nu_i\cup R_3\nu'_i\cup R_3\mu'_i$, $i=1,2$, in the
torus $\Cal T:=\partial\partial_0\Cal P=\partial(\Cal P\cap\S V)$.

The $4$-Cake is tessellated by $8$ conjugates to $\Cal P$ :
$$\Cal P_0:=W_0\Cal P,\qquad\Cal P_1:=W_2R_3\Cal P,\qquad\Cal
P_2:=W_4\Cal P,\qquad\Cal P_3:=W_5\Cal P,$$
$$\Cal P_4:=W_6\Cal P,\qquad\Cal P_5:=W_8R_3\Cal P,\qquad\Cal
P_6:=W_{10}\Cal P,\qquad\Cal P_7:=W_{11}\Cal P.$$

\vskip5pt

\noindent
$\hskip0pt\vcenter{\hbox{\epsfbox{Pict.4}}}$

\newpage

As is shown in Subsection I.2.4, $\partial_1\Cal C:=\Cal C\cap\S V$ and
$\partial_1\Cal P_i:=\Cal P_i\cap\S V$, $0\le i\le7$, are solid tori
and the fundamental group
$\pi_1\partial_1\Cal C=\pi_1\partial_1\Cal P_i$ is
generated\footnote{We denote by $[c]$ the homotopy class of a closed
curve $c$.}
by $g:=[C_1\cap\S V]$, where the circle $C_1\cap\S V$ is oriented with
respect to the natural orientation of $C_1\cap\overline\B V$. Since
$\pi_1\partial_1\Cal C\simeq\Bbb Z$ is abelian, we can write
$g=[C\cap\S V]$ for any slice $C$ of any bisector involved in one of
the $\Cal P_i$'s (we orient the circle $C\cap\S V$ naturally).
Therefore,
$$W_0g=g,\qquad W_2R_3g=g,\qquad W_4g=g,\qquad W_5g=g,$$
$$W_6g=g,\qquad W_8R_3g=g,\qquad W_{10}g=g,\qquad W_{11}g=g.$$

For suitable $n_i\in\Bbb Z$, we have $[\gamma_i]=n_ig$, $i=1,2$.

\medskip

{\bf3.2. Lemma.} {\sl The closed curve
$$\nu_1\cup R_3\nu'_1\cup W_2R_3\mu_2\cup W_2R_3\nu_2\cup
W_4R_3\mu'_1\cup W_4\mu_1\cup W_5R_3\nu'_2\cup W_5R_3\mu'_2\cup$$
$$\cup W_6\nu_1\cup W_6R_3\nu'_1\cup W_8R_3\mu_2\cup
W_8R_3\nu_2\cup W_{10}R_3\mu'_1\cup W_{10}\mu_1\cup
W_{11}R_3\nu'_2\cup W_{11}R_3\mu'_2\subset\partial\partial_0\Cal C$$

\smallskip

\noindent
provides a section of the disc bundle\/ $M$. The Euler number\/
$eM$ equals\/ $4n_1+4n_2$.}

\medskip

{\bf Proof} consists mostly in looking at the picture above Lemma 3.2,
where the $4$-Cake is cut into its pieces for the sake of convenience.
The reader is supposed to check that the $4$ conjugates to $\gamma_1$
and the $4$ conjugates to $\gamma_2$ located in the corresponding
conjugates to $\Cal T$ suit each other with respect to gluing back the
$4$-Cake and with respect to the isometries $I_i$'s. The rest follows
from the analog of Remark~I.2.4.4 adapted to our case
$_\blacksquare$

\medskip

We put $\zeta_i:=\mu_i\cup\nu_i\cup\ell\subset\partial\partial_0P$,
$i=1,2$.

\medskip

{\bf3.3 Lemma.} {\sl$[\gamma_i]=2[\zeta_i]$.}

\medskip

{\bf Proof.} Since $R_3$ alters the natural orientation of every
complex geodesic, we obtain $g+R_3g=0$. Hence,
$[\zeta_i]+[R_3\zeta_i]=0$, $i=1,2$. It is immediate that
$[\gamma_i]=[\zeta_i]+[R_3\zeta'_i]=2[\zeta_i]$
$_\blacksquare$

\medskip

We will specify the curves $\mu_1$ and $\nu_1$. We draw the ideal
meridional curve $\ell_1$ of $B[C_1,C_2]$ that begins with $w_1$ and
let $f_1\in C_2\cap\S V$ denote its end. The curve $\ell_1$ intersects
$M_1$ in some $k_1\in M_1\cap\S V$ and thus is divided into two parts,
$\ell_1=\ell_1^1\cup\ell_1^2$. We also draw the ideal meridional curve
$\ell_2$ of $B[C_2,C_3]$ that ends with $w_3$. Let $f_2$ denote the
beginning of $\ell_2$. Similarly, $\ell_2$ intersects $M_2$ in some
$R_2k_2\in M_2\cap\S V$ for a suitable $k_2\in M_2\cap\S V$ and
$\ell_2$ is divided into two parts $\ell_2=\ell_2^1\cup\ell_2^2$. We
denote by $h$ the arc in the circle $C_2\cap\S V$ that, following the
natural orientation of $C_2\cap\S V$, begins with $f_1$ and ends with
$f_2$. According to the definition at the beginning of Subsection
I.2.5, the curve $c:=\ell_2\cup\ell\cup\ell_1\cup h$ is standard since
$\Delta$ is an elliptic transversal triangle oriented in the
counterclockwise sense. By Theorem I.2.5.1, $c$~is a trivializing
curve, i.e., $[c]=0$.

By Lemma I.2.1.1, $f_1=R(q_1)R(q_3)w_3$. By the properties of the
triangle $\Delta$, the point $f_1$ lies on the side of the normal
vector to $\G{\prec}b_2,e_2{\succ}$. We draw the ideal meridional curve
$\ell_3$ of $B[M_1,C_2]$ that begins with $R_1k_1\in M_1\cap\S V$. The
reflection $R_1$ interchanges the ideal meridional curves $\ell_1^1$
and $\ell_3$. It follows from $w_2=R_1w_1$ that $\ell_3$ ends with
$w_2$. Similarly, the ideal meridional curve $\ell_4$ of $B[C_2,M_2]$
that begins with $w_2$ should end with $k_2$. We denote by $h_j^i$,
$i=1,2$, $j=1,2$, the arc of the circle $M_j\cap\S V$ that begins with
$k_j$ and ends with $R_jk_j$ such that $h_j^1$ follows the natural
orientation of the circle and $h_j^2$ goes in opposite direction.

We put $\mu_1:=\ell_1^1\cup h_1^1\cup\ell_3$ and
$\nu_1:=\ell_4\cup h_2^1\cup\ell_2^2$. It is easy to see that
$\mu_2:=R_1\mu'_1=\ell_1^1\cup h_1^2\cup\ell_3$ and
$\nu_2:=R_2\nu'_1=\ell_4\cup h_2^2\cup\ell_2^2$.

\noindent
$\hskip0pt\vcenter
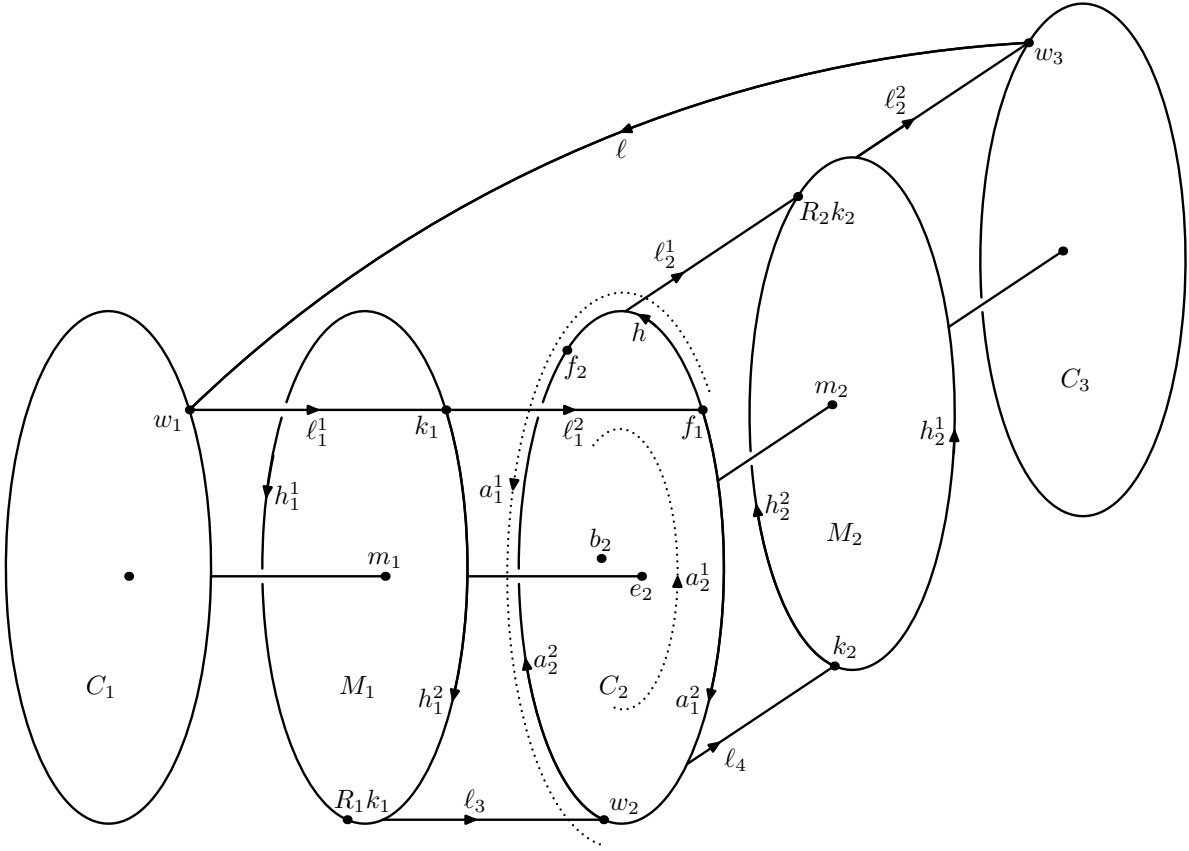
{\hbox{\epsfbox{Pict.5}}}$

\vskip10pt

We define 4 arcs of the circle $C_2\cap\S V$. The arcs $a_1^i$,
$i=1,2$, begin with $f_1$ and end with $w_2$. The arcs~$a_2^i$,
$i=1,2$, begin with $w_2$ and end with $f_2$. The arcs $a_j^1$ follow
the natural orientation of $C_2\cap\S V$, and the arcs $a_j^2$ do not.
It is easy to see that the curve $h_1^i\cup\ell_3$ is homotopic in the
cylinder $B[M_1,C_2]\cap\S V\subset\Cal T$ to $\ell_1^2\cup a_1^i$ and
that the curve $\ell_4\cup h_2^i$ is homotopic in the cylinder
$B[C_2,M_2]\cap\S V\subset\Cal T$ to $a_2^i\cup\ell_2^1$, $i=1,2$. This
implies that the curve
$$\zeta_i=\mu_i\cup\nu_i\cup\ell=\ell_1^1\cup
h_1^i\cup\ell_3\cup\ell_4\cup h_2^i\cup\ell_2^2\cup\ell$$

\noindent
$\hskip312pt\vcenter{\hbox{\epsfbox{Pict.6}}}$

\rightskip150pt

\vskip-144pt

\noindent
is homotopic in $\Cal T$ to the curve

\vskip12pt

\hskip23pt$\ell_1^1\cup\ell_1^2\cup a_1^i\cup
a_2^i\cup\ell_2^1\cup\ell_2^2\cup\ell'=\ell_1\cup a_1^i\cup
a_2^i\cup\ell_2\cup\ell$.

\vskip9pt

\noindent
From $0=[c]=[\ell_2\cup\ell\cup\ell_1\cup h]$, we conclude that
$[\zeta_i]=[a_1^i\cup a_2^i\cup h']$.

Clearly, $R(m_1)R(q_1)f_1=w_2$. Since $R(m_1)R(q_1)$ is the reflection
in the complex spine of $B[C_1,C_2]$ and $f_1,w_2\in C_2$, we obtain
$R(e_2)f_1=\nomathbreak w_2$. In a similar way, we can show that
$R(b_2)f_2=w_2$. Therefore, $f_1$ and $f_2$ lie on the same side from
$\G{\prec}b_2,e_2{\succ}$ and $h'$ is entirely included into the same
part of the circle $C_2\cap\S V$. This implies that
$[\zeta_1]=[C_2\cap\S V]$ and $[\zeta_2]=-[C_2\cap\S V]$. We arrived at
the

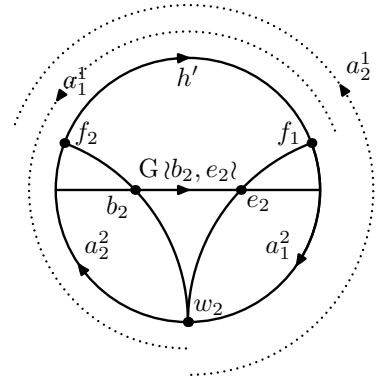
\rightskip0pt

\medskip

{\bf3.4. Theorem.} {\sl$eM=0$
$_\blacksquare$}

\medskip

{\bf3.5. Remark.} It is easy to see that we could have arrived at
$eM=-16$, would $f_1$ lie on the other side from
$\G{\prec}b_2,e_2{\succ}$.

\medskip

{\bf3.6. Remark.} In order to get a trivial complex hyperbolic disc
bundle over an orientable surface of genus $2$, we take the group $H_5$
generated by
$$X_1:=R_1,\qquad X_2:=R_2,\qquad X_3:=R_3R_2R_3,\qquad
X_4:=R_3R_1R_3,\qquad X_5:=R_0.$$
The defining relations of $H_5$ are $X_5X_4X_3X_2X_1=1$ and $X_i^2=1$.
The discreteness of $H_5$ is implied by that of $G$. Passing to a
suitable subgroup of index 4 in $H_5$, we arrive at a desired bundle
(see the plane example in Subsection I.2.1 with $n=5$).

\medskip

{\bf3.7. Remark.} A few comments about [Eli, Open Question 8.1]. By
Lemma I.2.2.2, the bisectors forming $\partial_0{\Cal C}$ remain
transversal along their common slices in some open neighbourhood $U$ of
$\overline\B V$. Also, they are transversal to $\S V$. It follows that
these bisectors limit in $U$ some {\it extended\/} $4$-Cake
$\hat\Cal C$ and that a finite number of conjugate to $\hat\Cal C$
tessellate some open neighbourhood of $\hat\Cal C$. With routine use of
the fact that $\Cal C\cap\S V$ is a compact (it is a solid torus), we
can diminish $U$ and reach a situation where the $I_i$'s identify the
parts of the corresponding bisectors that belong to $\hat\Cal C$. In
other words, we construct some holomorphic (open) manifold $N$ such
that $M\subset N$ and the ideal boundary of $M$, homeomorphic to
$\Sigma\times\Bbb S^1$, is simply $\overline M\setminus M$.

\medskip

As for the complex hyperbolic variant of the problem in [GLT], it seems
that cooking a cake is insufficient to solve it, or that the
ingredients are wrong.

\bigskip

\centerline{\bf4. Existence of the Triangle: A Piece of Complex
Hyperbolic Geometry}

\bigskip

In this section, we deal with geometric aspects related to the
existence of the triangle described in Subsection 2.2. The most
important part, constructing the meridian $\Gamma$ of the side
$B[C_3,C_1]$, explores some curious geometric ideas.

\medskip

{\bf4.1. Lemma.} {\sl Let\/ $C$ be a complex geodesic, let\/ $\G$ be a
geodesic not included in\/ $C$, and let\/ $g$ and\/ $g'$ be
different negative points in $\G$. Denote by\/ $p$ the polar point to\/
$C$. If\/ $\G\cap C\cap\B V\ne\varnothing$, then\/
$\displaystyle{\frac{\langle g,p\rangle\langle p,g'\rangle}{\langle
g,g'\rangle}}\in\Bbb R$.
Suppose that\/
$\displaystyle{\frac{\langle g,p\rangle\langle p,g'\rangle}{\langle
g,g'\rangle}}\notin\Bbb R$.
Then there exists a unique point in\/ $\G$ closest to\/ $C$ and the
condition that\/ $\ta(\G,C)=\ta(g,C)$ is equivalent to\/
$\Re\displaystyle{\frac{\langle p,g'\rangle\langle g,g\rangle}{\langle
g,g'\rangle\langle p,g\rangle}}=1$.}

\medskip

{\bf Proof.} Let $g''\in\G\cap C\cap\B V$ and let $q$ stand for the
polar point to the complex geodesic $L$ spanned by $\G$. It follows
from $\G\not\subset C$ that $p\ne q$. Therefore, we can write $p=cq+l$,
where $c\in\Bbb C$ and $l$ represents a point in $L$. The point $g''$
is orthogonal to $l$ since $g''$ is orthogonal to $p$ and to $q$. By
Remark~I.5.1.3, $l$~belongs to the extended geodesic related to $\G$.
Now, Lemma I.5.1.5 implies the first assertion.

Suppose that
$\displaystyle{\frac{\langle g,p\rangle\langle p,g'\rangle}{\langle
g,g'\rangle}}\notin\Bbb R$.
Denote by $v_1,v_2$ the vertices of $\G$. In view of $\G\not\subset C$,
we can assume that $v_1\notin C$ and take representatives
$p,v_1,v_2\in V$ with the Gram matrix
$\left(\smallmatrix 1&1&z\\1&0&-\frac12\\\overline
z&-\frac12&0\endsmallmatrix\right)$,
$z\in\Bbb C$. By~Lemma~I.5.1.6, every point in $\G\cap\B V$ has the
form
$$g(x):=xv_1+x^{-1}v_2,\quad x>0,\qquad\big\langle
g(x),g(x)\big\rangle=-1.$$
In particular, $g'=g(x')$ for a suitable $x'>0$. From
$\displaystyle{\frac{\langle g,p\rangle\langle p,g'\rangle}{\langle
g,g'\rangle}}\notin\Bbb R$,
we conclude that $z\ne0$. By~Lemma~I.5.2.10,
$$\ta\big(g(x),C\big)=1-\ta\big(g(x),p\big)=
1+(x+x^{-1}\overline z)(x+x^{-1}z)=1+2\Re z+x^2+|z|^2x^{-2}.$$
This function in $x$ takes its only minimum when $x=\sqrt{|z|}$. We
have
$$\langle p,g'\rangle\big\langle
g(x),g(x)\big\rangle=-(x'+{x'}^{-1}z),\qquad\big\langle
g(x),g'\big\rangle\big\langle
p,g(x)\big\rangle=-\frac{(x{x'}^{-1}+x^{-1}x')(x+x^{-1}z)}{2}.$$
Therefore,
$$\Re\frac{\langle p,g'\rangle\big\langle
g(x),g(x)\big\rangle}{\big\langle g(x),g'\big\rangle\big\langle
p,g(x)\big\rangle}=\Re\frac{2(x'+{x'}^{-1}z)}
{(x{x'}^{-1}+x^{-1}x')(x+x^{-1}z)}=$$
$$=\frac{(x{x'}^{-1}+x^{-1}x')(x^2+x^{-2}|z|^2+2\Re
z)+x^{-3}{x'}^{-1}(x^2-{x'}^2)\big(|z|^2-x^4\big)}
{(x{x'}^{-1}+x^{-1}x')(x^2+x^{-2}|z|^2+2\Re z)}=$$
$$=1+\frac{x^{-3}{x'}^{-1}(x^2-{x'}^2)\big(|z|^2-x^4\big)}
{(x{x'}^{-1}+x^{-1}x')(x^2+x^{-2}|z|^2+2\Re z)}.$$
Consequently,
$\Re\displaystyle{\frac{\langle p,g'\rangle\big\langle
g(x),g(x)\big\rangle}{\big\langle g(x),g'\big\rangle\big\langle
p,g(x)\big\rangle}}=1$
if and only if $(x^2-{x'}^2)\big(|z|^2-x^4\big)=0$
$_\blacksquare$

\medskip

{\bf4.2. Lemma {\rm(folklore)}.} {\sl Suppose that\/ $\tr I>3$ for
some\/ $I\in\SU V$. Then there exist a geodesic\/ $\G$ and two distinct
points\/ $g,g'\in\G\cap\B V$ such that\/ $I=R(g')R(g)$. Moreover, we
can arbitrarily choose\/ $g\in\G\cap\B V$.}

\medskip

{\bf Proof.} The $\Bbb C$-vector space $V$ admits a basis $p,v_1,v_2$
with the Gram matrix
$\left(\smallmatrix1&0&0\\0&0&-\frac12\\0&-\frac12&0\endsmallmatrix
\right)$.
There exists a unique $r>1$ such that $\tr I=1+r+r^{-1}$. The matrix
$\left(\smallmatrix1&0&0\\0&r&0\\0&0&r^{-1}\endsmallmatrix\right)$
defines some $J\in\SU V$. By [Gol2, Theorem 6.2.4, p.~204], the
isometries $I$ and $J$ are conjugate. We put $\G:=\G[v_1,v_2]$. Let
$g=\alpha v_1+\alpha^{-1}v_2$, $\alpha>0$, be an arbitrary point in
$\G\cap\B V$. A straightforward verification shows that $J=R(g')R(g)$,
where $g':=\alpha\sqrt rv_1+(\alpha\sqrt r)^{-1}v_2$
$_\blacksquare$

\medskip

We postpone until Section 5 the proof of the following three algebraic
lemmas.

\medskip

{\bf4.3. Lemma {\rm(compare with [Pra])}.} {\sl Let\/
$x_1,x_2,x_3\in V$ be nonisotropic. Then
$$\big\langle R(x_2)x_1,x_1\big\rangle=\big(2\ta(x_1,x_2)-1\big)\langle
x_1,x_1\rangle,\qquad\tr\big(R(x_2)R(x_1)\big)=4\ta(x_1,x_2)-1,$$
$$\tr\big(R(x_3)R(x_2)R(x_1)\big)=8\frac{\langle x_1,x_2\rangle\langle
x_2,x_3\rangle\langle x_3,x_1\rangle}{\langle x_1,x_1\rangle\langle
x_2,x_2\rangle\langle
x_3,x_3\rangle}-4\ta(x_1,x_2)-4\ta(x_2,m_3)-4\ta(x_3,x_1)+3.$$}

{\bf4.4. Lemma.} {\sl Suppose that the Gram matrix of\/
$p_1,p_2,p_3\in V$ has the form\/
$\left(\smallmatrix1&t_1&t\\t_1&1&t_2\overline\lambda\\t&t_2\lambda&1
\endsmallmatrix\right)$,
where\/ $t,t_1,t_2>1$, $|\lambda|=1$, and\/ $\lambda\ne1$. We put\/
$m_1:=\displaystyle\frac{p_1-p_2}{\sqrt{2(t_1-1)}}$ and\/
$m_2:=\displaystyle\frac{\lambda p_2-p_3}{\sqrt{2(t_2-1)}}$. Then
$$\Re\frac{\langle p_1,m_2\big\rangle\langle m_1,m_1\rangle}{\langle
m_1,m_2\rangle\langle
p_1,m_1\rangle}=1+\frac{t^2-tt_1+t_1^2-(t_2-1)^2+tt_1(1-2\Re\lambda)}
{(t_1+t_2-1)^2+t^2-2t(t_1+t_2-1)\Re\lambda},$$
$$\tr\big(R(m_2)R(m_1)R(p_1)\big)=2t(\overline\lambda-1)+
\frac{2tt_1-t-t_1+1-2t_2}{t_1-1}-$$
$$-\frac{t^2-tt_1+t_1^2-(t_2-1)^2+t(t_1+t_2-1)
(1-2\Re\lambda)}{(t_1-1)(t_2-1)}.$$}

\medskip

{\bf4.5. Lemma.} {\sl For every\/ $t>\frac32$, there exist unique\/
$t_1>1$ and\/ $t_2\in\Bbb R$ such that
$$t^2-tt_1+t_1^2-(t_2-1)^2=0,\qquad2tt_1-t-t_1+1-2t_2=0.\eqno{(1)}$$
For such\/ $t,t_1,t_2$, we have\/ $t_2>t_1$.}

\medskip

{\bf4.6. Proposition.} {\sl Suppose that the Gram matrix of\/
$p_1,p_2,p_3\in V$ has the form\/
$\left(\smallmatrix1&t_1&t\\t_1&1&t_2\overline\vartheta\\
t&t_2\vartheta&1\endsmallmatrix\right)$,
where\/ $\vartheta:=\exp\displaystyle\frac{\pi i}3$ and\/
$t,t_1,t_2\in\Bbb R$ satisfy the equalities\/ {\rm(1)} and the
inequalities\/ $t>\frac32$, $t_1>1$. We put
$$C_i:=\Bbb{CP}p_i^\perp\cap\overline\B V,\quad i=1,2,3,\qquad
m_1:=\frac{p_1-p_2}{\sqrt{2(t_1-1)}},\qquad m_2:=\frac{\vartheta
p_2-p_3}{\sqrt{2(t_2-1)}},$$
$$R_i:=R(m_i),\quad i=1,2,\qquad R_0:=R(p_1).$$
Then\/ $m_i$ is the middle point of the real spine of the segment\/
$B[C_i,C_{i+1}]$, $i=1,2$, and there exists a meridian\/ $\Gamma$ of
the segment\/ $B[C_3,C_1]$ such that\/
$R_3R_1R_2R_3R_2R_1R_0=\vartheta^{-2}$, where\/ $R_3$ denotes the
reflection in $\Gamma$.}

\medskip

{\bf Proof.} The fact that $m_i$ is the middle point of the real spine
of the segment $B[C_i,C_{i+1}]$ is straightforward.

Since $2\Re\vartheta=1$ and $\vartheta^2(\overline\vartheta-1)=1$, we
obtain $\tr(\vartheta^2R_2R_1R_0)=2t>3$ by Lemma 4.4. It follows from
$\vartheta^2\in\SU V$ that $\vartheta^2R_2R_1R_0\in\SU V$. By Lemma
4.2, $\vartheta^2R_2R_1R_0=R(m'_2)R(m'_1)$ for suitable distinct
$m'_1,m'_2\in\B V$.

It follows from $m_1,m_2\in\B V$ that $\ta(m_1,m_2)\ge1$. Suppose that
$\displaystyle{\frac{\langle m'_1,p_1\rangle\langle
p_1,m'_2\rangle}{\langle m'_1,m'_2\rangle}}\in\Bbb R$.
Then, by~Lemma 4.3, $0\ne\tr(R_2R_1)\in\Bbb R$ and
$\tr\big(R(m'_2)R(m'_1)R_0\big)\in\Bbb R$. This contradicts the
equality $\vartheta^2R_2R_1=R(m'_2)R(m'_1)R_0$. Hence,
$\displaystyle{\frac{\langle m'_1,p_1\rangle\langle
p_1,m'_2\rangle}{\langle m'_1,m'_2\rangle}}\notin\Bbb R$.
In particular, $m'_1,m'_2\notin C_1$. By Lemmas 4.1 and 4.2, we~can
assume that
$\Re\displaystyle\frac{\langle p_1,m'_2\big\rangle\langle
m'_1,m'_1\rangle}{\langle m'_1,m'_2\rangle\langle p_1,m'_1\rangle}=1$.

We denote $p'_2:=-R(m'_1)p_1\in V$ and $C'_2:=\Bbb{CP}{p'_2}^\perp$. If
$\langle p'_2,m'_2\rangle=0$, then
$\big\langle p_1,R(m'_1)m'_2\big\rangle=0$ and
$R(m'_1)m'_2\in\G{\prec}m'_1,m'_2{\succ}\cap C_1\cap\B V$. By Lemma
4.1, this is impossible since
$\displaystyle{\frac{\langle m'_1,p_1\rangle\langle
p_1,m'_2\rangle}{\langle m'_1,m'_2\rangle}}\notin\Bbb R$.
Therefore, by~Lemma 4.3,
$$t'_1:=\langle p'_2,p_1\rangle=-2\ta(p_1,m'_1)+1>1,\qquad
t'_2:=-\big\langle R(m'_2)p'_2,p'_2\big\rangle=-2\ta(p'_2,m'_2)+1>1.$$
It is easy to verify that $\langle m_i,m_i\rangle=-1$ for $i=1,2$,
that $R_1p_1=-p_2$, and that $R_2\vartheta p_2=-p_3$. Hence, the
equalities $R_0p_1=p_1$ and $\vartheta^2R_2R_1R_0=R(m'_2)R(m'_1)$
imply $\vartheta p_3=-R(m'_2)p'_2$. In other words, the Gram matrix of
$p_1,p'_2,p_3$ has the form
$\left(\smallmatrix1&t'_1&t\\t'_1&1&t'_2\vartheta\\
t&t'_2\overline\vartheta&1\endsmallmatrix\right)$
with $t'_1,t'_2>1$.

It follows from $p'_2:=-R(m'_1)p_1$ and $\vartheta p_3=-R(m'_2)p'_2$
that $R(m'_1)C_1=C'_2$ and $R(m'_2)C'_2=C_3$. Consequently, $m'_1$ and
$m'_2$ are the middle points of the real spines of the segments
$B[C_1,C'_2]$ and $B[C'_2,C_3]$. So, we can assume that
$$m'_1=\frac{p_1-p'_2}{\sqrt{2(t'_1-1)}},\qquad
m'_2=\frac{\overline\vartheta p'_2-p_3}{\sqrt{2(t'_2-1)}}.$$
By Lemma 4.4, we conclude from
$\Re\displaystyle\frac{\langle p_1,m'_2\big\rangle\langle
m'_1,m'_1\rangle}{\langle m'_1,m'_2\rangle\langle p_1,m'_1\rangle}=1$
and from $2\Re\overline\vartheta=1$ that
$$t^2-tt'_1+{t'_1}^2-(t'_2-1)^2=0.$$

Taking into account that $2\Re\vartheta=1$, we obtain
$$\langle m_1,m_2\rangle\langle
m_2,m_1\rangle=\frac{\big(\overline\vartheta(t_1+t_2-1)-
t\big)\big(\vartheta(t_1+t_2-1)-t\big)}{4(t_1-1)(t_2-1)}=$$
$$=\frac{(t_1+t_2-1)^2+t^2-2t(t_1+t_2-1)\Re\vartheta}
{4(t_1-1)(t_2-1)}=$$
$$=\frac{2t+1}4+\frac{t^2-tt_1+t_1^2-(t_2-1)^2-(2tt_1-t-t_1+1-2t_2)
(t_2-1)}{4(t_1-1)(t_2-1)}.$$
The equalities (1) imply that $\ta(m_1,m_2)=\displaystyle\frac{2t+1}4$.
By Lemma 4.3, $\tr(R_2R_1)=2t$. From
$\vartheta^2R_2R_1=R(m'_2)R(m'_1)R_0$, we derive
$\tr\big(R(m'_2)R(m'_1)R_0\big)=2t\vartheta^2=2t(\vartheta-1)$. By
Lemma 4.4,
$$2tt'_1-t-t'_1+1-2t'_2=0.$$
By Lemma 4.5, $t'_1=t_1$ and $t'_2=t_2$.

The Gram matrices of the triples $p_1,p_2,p_3$ and $p_1,p'_2,p_3$ are
complex conjugate. Therefore, there exists a real plane
$\Gamma=\Bbb{CP}S$ such that $p_1,p_3\in S$ and $R_3p_2=p'_2$, where
$S\subset V$ is a three-dimensional real subspace with hermitian form
real over it, $V=S\otimes_\Bbb R\Bbb C$, and $R_3\in\widehat\SU V$
stands for the reflection in~$S$ (induced by the complex conjugation in
$\Bbb C$). Hence, $R_3m_i=m'_i$ and, therefore, $R(m'_i)=R_3R_iR_3$,
$i=1,2$. The relation $\vartheta^2R_2R_1R_0=R(m'_2)R(m'_1)$ implies
the relation $R_3R_1R_2R_3R_2R_1R_0=\vartheta^{-2}$
$_\blacksquare$

\medskip

{\bf4.7. Corollary.} {\sl In terms of Proposition\/ {\rm4.6}, let\/
$\G[d_3,c_1]$ stand for the real spine of the segment\/
$B[C_3,C_1]$, where\/ $c_1\in C_1$ and\/ $d_3\in C_3$. We put
$$c_2:=R_1c_1,\quad c_3:=R_2c_2,\quad d_2:=R_2d_3,\quad
d_1:=R_1d_2,\quad S_i:=\G[c_i,c_{i+1}],\quad T_i:=\G[d_i,d_{i+1}],$$
$i=1,2,3$ {\rm(}the indices are modulo $3${\rm).} Suppose that\/
$p_1,p_2,p_3$ form a\/ $\Bbb C$-basis of\/ $V$. Then
$$\#S_i\cap T_i\le1,\quad i=1,2,3,\qquad
c_3,d_1\in\Gamma,\qquad\frac{\langle c_1,c_3\rangle\langle
c_3,c_2\rangle}{\langle c_1,c_2\rangle}\in\Bbb R\overline\vartheta i.$$
If the triangle\/ $\Delta:=\Delta(C_1,C_2,C_3)$ is transversal, then it
is oriented in the counterclockwise sense and the triangles\/ $\Delta$
and\/ $R_3\Delta$ lie on different sides from\/
$B{\prec}C_3,C_1{\succ}$.}

\medskip

{\bf Proof.} Suppose that $\#S_i\cap T_i\ge2$ for some $i$. Then
$S_i=T_i$, which implies that $c_j=d_j$ for all $j$. In particular,
$d_3=R_2R_1c_1$. We can take $c_1=p_3-tp_1$. Obviously,
$d_3\in\Bbb Cp_1+\Bbb Cp_3$. Hence, the coefficient of $p_2$ in
$R_2R_1c_1$ vanishes. By a straightforward calculus, this coefficient
equals
$$\frac{tt_1-t^2\vartheta-t_1t_2\vartheta+tt_2\vartheta^2}
{(t_1-1)(t_2-1)}+\vartheta.$$
Taking into account that $\vartheta^2=\vartheta-1$, we obtain
$$\frac{t(t_1-t_2)}
{(t_1-1)(t_2-1)}+\Big(\frac{tt_2-t^2-t_1t_2}
{(t_1-1)(t_2-1)}+1\Big)\vartheta=0$$
and conclude that $t_1=t_2$. This contradicts Lemma 4.5.

We take representatives $c_1,d_3\in S$. It follows from $R_0c_1=-c_1$,
$R_3R_2R_1R_0=\vartheta^{-2}R_2R_1R_3$, and $R_3c_1=c_1$ that
$$R_3c_3=R_3R_2R_1c_1=-R_3R_2R_1R_0c_1=-\vartheta^{-2}R_2R_1R_3c_1=
-\vartheta^{-2}R_2R_1c_1=-\vartheta^{-2}c_3.$$
Similarly, from $R_0R_3d_1=-R_3d_1$, $R_3d_3=d_3$, and
$R_0R_3R_1R_2R_3=\vartheta^{-2}R_1R_2$, we obtain
$$R_3d_1=-R_0R_3d_1=-R_0R_3R_1R_2d_3=-R_0R_3R_1R_2R_3d_3=
-\vartheta^{-2}R_1R_2d_3=-\vartheta^{-2}d_1.$$
Therefore, $c_3,d_1\in\Gamma$.

It is immediate that $R_3(cx)=\overline c R_3x$ for all $c\in\Bbb C$
and $x\in V$. Hence, the equality $R_3c_3=-\vartheta^{-2}c_3$ implies
the equality $R_3(\overline\vartheta ic_3)=\overline\vartheta ic_3$. In
other words, $\overline\vartheta ic_3\in S$. Consequently,
$\langle c_1,c_3\rangle\in\Bbb R\overline\vartheta i$. From
$$\langle c_1,c_2\rangle=\big\langle c_1,R_1c_1\rangle=\big\langle
c_1,-2\langle c_1,m_1\rangle m_1-c_1\big\rangle=-2\langle
c_1,m_1\rangle\langle m_1,c_1\rangle-\langle c_1,c_1\rangle\in\Bbb R,$$
$$\langle c_3,c_2\rangle=\big\langle R_3c_2,c_2\rangle=\big\langle
-2\langle c_2,m_2\rangle m_2-c_2,c_2\big\rangle=-2\langle
c_2,m_2\rangle\langle m_2,c_2\rangle-\langle c_2,c_2\rangle\in\Bbb R,$$
we conclude that
$\displaystyle\frac{\langle c_1,c_3\rangle\langle
c_3,c_2\rangle}{\langle c_1,c_2\rangle}\in\Bbb R\overline\vartheta i$.

Suppose that $\Delta$ is transversal. By Criterion I.2.5.4, the
triangles $\Delta$ and $R_3\Delta'$ are oriented in the
counterclockwise sense since $\Im\overline\vartheta<0$
$_\blacksquare$

\bigskip

\centerline{\bf5. Existence of the Triangle: A Piece of Elementary
Algebra}

\bigskip

In this section, we prove Lemmas 4.3--5 and complete the proof of the
existence of the triangle.

\medskip

{\bf Proof of Lemma 4.3.} The equality
$\big\langle R(x_2)x_1,x_1\big\rangle=\big(2\ta(x_1,x_2)-1\big)\langle
x_1,x_1\rangle$
is straightforward.

As is easy to see,
$\tr\big(\langle -,x\rangle y\big)=\langle y,x\rangle$ for the
$\Bbb C$-linear transformation $\langle -,x\rangle y:V\to V$ given by
the rule $\langle -,x\rangle y:v\mapsto\langle v,x\rangle y$. Denote
$\varphi_i:=2\displaystyle\frac{\langle -,x_i\rangle x_i}{\langle
x_i,x_i\rangle}$.
Then
$$\tr1=3,\qquad\tr\varphi_i=2,\qquad\tr(\varphi_j\varphi_i)=4\ta(x_i,x_j),\qquad
\tr(\varphi_k\varphi_j\varphi_i)=8\frac{\langle x_i,x_j\rangle\langle
x_j,x_k\rangle\langle x_k,x_i\rangle}{\langle x_i,x_i\rangle\langle
x_j,x_j\rangle\langle
x_k,x_k\rangle}.$$
It remains to take into account that $R(x_i)=\varphi_i-1$
$_\blacksquare$

\medskip

{\bf Proof of Lemma 4.4} is straightforward. It is immediate that
$\langle m_i,m_i\rangle=-1$, $i=1,2$. Also, we have
$$\langle p_1,m_1\rangle=\frac{1-t_1}{\sqrt{2(t_1-1)}},\qquad\langle
p_1,m_2\rangle=\frac{\overline\lambda
t_1-t}{\sqrt{2(t_2-1)}},\qquad\langle
m_1,m_2\rangle=\frac{\overline\lambda(t_1+t_2-1)-t}{2\sqrt{(t_1-1)
(t_2-1)}}.$$
Clearly, $\lambda\ne1$ implies $\overline\lambda(t_1+t_2-1)-t\ne0$.
Hence,
$$\Re\frac{\langle p_1,m_2\big\rangle\langle m_1,m_1\rangle}{\langle
m_1,m_2\rangle\langle p_1,m_1\rangle}=-\Re\frac{2(\overline\lambda
t_1-t)(t_1-1)}{\big(\overline\lambda(t_1+t_2-1)-t\big)(1-t_1)}=
\Re\frac{2(t_1-\lambda t)}{t_1+t_2-1-\lambda t}=$$
$$=\Re\frac{2(t_1-\lambda t)(t_1+t_2-1-\overline\lambda
t)}{(t_1+t_2-1)^2+t^2-t(t_1+t_2-1)(\lambda+\overline\lambda)}=
\frac{2t_1(t_1+t_2-1)+2t^2-2t(2t_1+t_2-1)\Re\lambda}
{(t_1+t_2-1)^2+t^2-2t(t_1+t_2-1)\Re\lambda}=$$
$$=1+\frac{(t_1-t_2+1)(t_1+t_2-1)+t^2-2tt_1\Re\lambda}{(t_1+t_2-1)^2+
t^2-2t(t_1+t_2-1)\Re\lambda}=1+\frac{t^2-tt_1+t_1^2-(t_2-1)^2+tt_1
(1-2\Re\lambda)}{(t_1+t_2-1)^2+t^2-2t(t_1+t_2-1)\Re\lambda}.$$

By Lemma 4.3,
$$\tr\big(R(m_2)R(m_1)R(p_1)\big)=$$
$$=8\langle p_1,m_1\rangle\langle
m_1,m_2\rangle\langle m_2,p_1\rangle+4\langle p_1,m_1\rangle\langle
m_1,p_1\rangle-4\langle m_1,m_2\rangle\langle
m_2,m_1\rangle+4\langle m_2,p_1\rangle\langle p_1,m_2\rangle+3.$$
Therefore,
$$(t_1-1)(t_2-1)\tr\big(R(m_2)R(m_1)R(p_1)\big)=$$
$$=2(1-t_1)\big(\overline\lambda(t_1+t_2-1)-t\big)(\lambda t_1-t)+
2(t_2-1)(1-t_1)^2-$$
$$-\big(\overline\lambda(t_1+t_2-1)-t\big)\big(\lambda(t_1+t_2-1)-t\big
)+2(t_1-1)(\lambda t_1-t)(\overline\lambda t_1-t)+3(t_1-1)(t_2-1)=$$
$$=2(t_1-1)(\lambda t_1-t)\big(\overline\lambda(1-t_1-t_2)+t+
\overline\lambda t_1-t\big)+$$
$$+(t_1-1)(t_2-1)(2t_1-2+3)-(t_1+t_2-1)^2-t^2+
2t(t_1+t_2-1)\Re\lambda=$$
$$=2(t_1-1)(\lambda t_1-t)\overline\lambda(1-t_2)+(t_1-1)(t_2-1)
(2t_1+1)-(t_1+t_2-1)^2-\big(t^2-2t(t_1+t_2-1)\Re\lambda\big)=$$
$$=(t_1-1)(t_2-1)(2t\overline\lambda+1)-(t_1+t_2-1)^2-
\big(t^2-tt_1-tt_2+t+t(t_1+t_2-1)(1-2\Re\lambda)\big)=$$
$$=2t(t_1-1)(t_2-1)(\overline\lambda-1)+(t_1-1)(t_2-1)(2t+1)-
(t_1+t_2-1)^2+tt_2-t+t_1^2-(t_2-1)^2-$$
$$-\big(t^2-tt_1+t_1^2-(t_2-1)^2+t(t_1+t_2-1)(1-2\Re\lambda)\big)=$$
$$=2t(t_1-1)(t_2-1)(\overline\lambda-1)+(t_2-1)(2tt_1-t-t_1+1-2t_2)-$$
$$-\big(t^2-tt_1+t_1^2-(t_2-1)^2+t(t_1+t_2-1)(1-2\Re\lambda)\big)\
_\blacksquare$$

{\bf Proof of Lemma 4.5} is highly elementary. Using the second
equation in (1) for excluding $t_2$ from the first one, we arrive at
the equation $f(t_1)=0$, where
$$f(x):=(4t^2-4t-3)x^2-2(2t^2-t-1)x-(3t^2-2t-1)=$$
$$=(2t+1)(2t-3)x^2-2(2t+1)(t-1)x-(3t+1)(t-1).$$
Since $f(1)=-3t^2<0$ and $t>\frac32$, there is a unique solution
$t_1>1$. Explicitly,
$$t_1:=\frac{t-1}{2t-3}+
\frac2{2t-3}\sqrt{\frac{(2t^2-2t-1)(t-1)}{2t+1}}.\eqno{(2)}$$

Notice that
$$t_1=\frac{t-1}{2t-3}+\frac2{2t-3}\sqrt{\frac{(2t^2-2t-1)(t-1)}{2t+1}}
<\frac{t-1}{2t-3}+\frac2{2t-3}\sqrt{\frac{(2t^2-t-1)(t-1)}{2t+1}}=
\frac{3(t-1)}{2t-3}.$$

If the inequality $t_2>t_1$ is not true, we have
$2tt_1-t_1-t+1\le2t_1$. Therefore,
$t_1\le\displaystyle\frac{t-1}{2t-3}$. A~contradiction
$_\blacksquare$

\medskip

We fix some $t>\frac32$, define $t_1$ with (2), and put
$t_2:=\frac12(2tt_1-t-t_1+1)$. By Sylvester's Criterion, there exists
some $\Bbb C$-basis $p_1,p_2,p_3\in V$ with the Gram matrix
$G:=\left(\smallmatrix1&t_1&t\\t_1&1&t_2\overline\vartheta\\
t&t_2\vartheta&1\endsmallmatrix\right)$
if $\det G<0$. It follows from $\vartheta+\overline\vartheta=1$ that
$\det G=1-t^2-t_1^2-t_2^2+tt_1t_2$. Thus, the condition
$$t^2+t_1^2+t_2^2-tt_1t_2>1\eqno{(3)}$$
guarantees the existence of such a $\Bbb C$-basis $p_1,p_2,p_3\in V$.

We define the complex geodesics $C_i$, $i=1,2,3$, as in Proposition
4.6. Taking into account that $t_2>t_1$, we apply Criterion I.2.5.4 to
the triangle $\Delta:=\Delta(C_1,C_2,C_3)$. It follows from
$\Re\vartheta=\frac12$ that $\Delta$ is transversal if and only if
$$4tt_1t_2-t^2-4t_1^2-4t_2^2+4>0,\qquad
4tt_1t_2-4t^2-t_1^2-4t_2^2+4>0.\eqno{(4)}$$
By Corollary 4.7, $\Delta$ is oriented in the counterclockwise sense.

By Lemma I.5.4.3, $\Delta$ is elliptic if
$$\sqrt3\Big(1+\frac{t^2+t_1^2+t_2^2-tt_1t_2-1}
{(t+1)(t_1+1)(t_2+1)}\Big)<2.\eqno{(5)}$$

We take $m_i$ and $R_i$, $i=1,2$, as in Proposition 4.6 and define
$$c_1:=e_1:=p_3-tp_1,\qquad b_2:=p_3-t_2\vartheta p_2,\qquad
e_2:=p_1-t_1p_2,\qquad d_3:=b_3:=p_1-tp_3,$$
$$m_3:=p_3-p_1,\qquad q_1:=p_1+p_2,\qquad q_3:=p_3+p_1,\qquad
c_2:=R_1c_1,\qquad c_3:=R_2c_2.$$
It is immediate that these points are those described in Subsection
2.2.

We put
$$u:=\ta(c_3,d_3),\qquad w_3:=(u+\sqrt{u^2-u})c_3-\frac{\langle
c_3,d_3\rangle}{\langle d_3,d_3\rangle}d_3.$$
A straightforward verification shows that $w_3$ is isotropic.

Requiring that
$$u>1,\qquad\big\langle R(m_3)w_3,R_1R_2w_3\big\rangle\ne0,\qquad
\Im\frac{\big\langle b_2,R(q_1)R(q_3)w_3\big\rangle\big\langle
R(q_1)R(q_3)w_3,e_2\big\rangle}{\langle b_2,e_2\rangle}>0,\eqno{(6)}$$
we conclude that $w_3\ne0$, that $R(m_3)w_3\ne R_1R_2w_3$ (since both
$R(m_3)w_3$ and $R_1R_2w_3$ are isotropic), and that $b_2\ne e_2$.
Moreover, by Lemma I.5.2.9, the last of the inequalities (6) means that
the point $R(q_1)R(q_3)w_3$ lies in $C_2$ on the side of the normal
vector to the oriented geodesic $\G{\prec}b_2,e_2{\succ}$.

By Lemma I.3.1.3, the angle $\beta_i$ described in Subection 2.2 has
the form
$$\beta_i=\Arg\frac{\langle p_{i+1},c_i\rangle\langle
c_i,p_{i-1}\rangle}{\langle p_{i+1},p_i\rangle\langle
p_i,p_{i-1}\rangle},$$
$i=1,2,3$ (the indices are modulo $3$). Therefore,
$$\beta_1=\Arg\big(\langle p_2,c_1\rangle\langle
c_1,p_3\rangle\big),\qquad\beta_2=\Arg\big(\overline\vartheta\langle
p_3,c_2\rangle\langle c_2,p_1\rangle\big),\qquad
\beta_3=\Arg\big(\overline\vartheta\langle p_1,c_3\rangle\langle
c_3,p_2\rangle\big).$$
The fact that $\Delta$ is transversal and oriented in the
counterclockwise sense implies that $0<\beta_i<\pi$. Consequently, the
inequalities
$$\Re\big(\langle p_2,c_1\rangle\langle
c_1,p_3\rangle\big)>0,\qquad\Re\big(\overline\vartheta\langle
p_3,c_2\rangle\langle c_2,p_1\rangle\big)>0,\qquad
\Re\big(\overline\vartheta\langle p_1,c_3\rangle\langle
c_3,p_2\rangle\big)>0\eqno{(7)}$$
imply that $0<\beta_i<\frac\pi2$, $i=1,2,3$. Requiring in addition that
$$\Re\big(\overline\vartheta \langle p_2,c_1\rangle\langle
c_1,p_3\rangle\langle p_3,c_2\rangle\langle
c_2,p_1\rangle\big)>0,\eqno{(8)}$$
we obtain $\beta_1+\beta_2<\frac\pi2$, which implies
$0<\beta_1+\beta_2+\beta_3<\pi$.

Summarizing all the above, we arrive at the

\medskip

{\bf5.1. Proposition.} {\sl Under the conditions\/ {\rm(3--8),} the
triangle\/ $\Delta$ satisfies the properties listed in Subsection\/
{\rm2.2}\/
$_\blacksquare$}

\medskip

As we saw, every number involved in the inequalities (3--8) can be
explicitly expressed in terms of $t$.

\medskip

{\bf5.2. Lemma.} {\sl There exists\/ $t>\frac32$ satisfying the
inequalities {\rm(3--8).}}

\medskip

{\bf Proof} is absent. Using explicit expressions of everything
involved in the inequalities (3--8), the reader can check these
inequalities for an arbitrary $t\in[2.13,2.34]$. It suffices for the
needs of this article. For $t=2.22$, our calculator produced
$t_1=2.23$, $t_2=3.22$, and the following values on the left hand sides
of the inequalities:

\smallskip

\centerline{(3) 4.33; (4) 1.44, 1.56; (5) 1.86; (6) 1.43,
$-7.63-i\cdot4.41$, 3.68; (7) 13.11, 29.62, 31.05; (8) 248.24
$_\blacksquare$}

\medskip

Our attempt at working in terms of $\Gamma$ showed that the calculus
becomes much more bulky. It seems that the presented approach is more
economic.

\bigskip

\centerline{\bf6. References}

\bigskip

[AGG] S.~Anan$'$in, Carlos H.~Grossi, N.~Gusevskii, {\it Complex
Hyperbolic Structures on Disc Bundles over Surfaces. {\rm I.} General
Settings. A Series of Examples,} in preparation

[AGr] S.~Anan$'$in, Carlos H.~Grossi, {\it Local version of
Poincar\'e's polyhedron theorem,} in preparation

[Eli] Y.~Eliashberg, {\it Contact $3$-manifolds twenty years since
J.~Martinet's work,} Ann.~Inst.~Fourier (Grenoble) {\bf42} (1992),
No.~1--2, 165-192

[Gol1] W.~M.~Goldman, {\it Conformally flat manifolds with nilpotent
holonomy and the uniformization problem for $3$-manifolds,}
Trans.~Amer.~Math.~Soc.~{\bf278} (1983), No.~2, 573--583

[Gol2] W.~M.~Goldman, {\it Complex Hyperbolic Geometry,} Oxford
Mathematical Monographs. Oxford Science Publications. The Clarendon
Press, Oxford University Press, New York, 1999, xx+316 pp

[GLT] M.~Gromov, H.~B.~Lawson Jr., W.~Thurston, {\it Hyperbolic
$4$-manifolds and conformally flat $3$-manifolds,} Inst.~Hautes
\'Etudes Sci.~Publ.~Math., No.~68 (1988), 27--45

[Pra] A.~Pratoussevitch, {\it Traces in complex hyperbolic triangle
groups,} Geometriae Dedicata {\bf111} (2005), 159--185

[Sch] R.~E.~Schwartz, {\it Spherical {\rm CR} Geometry and Dehn
Surgery,} Ann.~of Math.~Stud. Princeton Univ. Press, Princeton, NJ, to
appear, available at http://www.math.brown.edu/$\sim$res/papers.html

[Tol] D.~Toledo, {\it Representations of surface groups in complex
hyperbolic space,} J.~Differential Geom. {\bf29} (1989), No.~1,
125--133

\enddocument